\documentclass[oneside,english,leqno]{amsart}
\usepackage[T1]{fontenc}
\usepackage[latin1]{inputenc}
\usepackage{amssymb}

\makeatletter
 \theoremstyle{plain}
\newtheorem{thm}{Theorem}[section]
  \theoremstyle{remark}
  \newtheorem{rem}[thm]{Remark}

\usepackage{babel}
\makeatother
\begin{document}

\title{Remarks on Ramanujan Function $A_{q}(z)$}

\author{Ruiming Zhang}

\curraddr{School of Mathematics\\
Guangxi Normal University\\
Guilin City, Guangxi 541004\\
P. R. China.}

\email{ruimingzhang@yahoo.com}

\begin{abstract}
In this short notes we will derive an inequality for scaled $q^{-1}$-Hermite
orthogonal polynomials of Ismail and Masson, an inequality for scaled
Stieltjes-Wigert, two inequalities for Ramanujan function and two
definite integrals for Ramanujan function. 
\end{abstract}



\maketitle

\section{Introduction\label{sec:Introduction}}

Ramanujan function $A_{q}(z)$, which is also called $q$-Airy function
in the literature, appears repeatedly in Ramanujan's work starting
from the Rogers-Ramanujan identities, where $A_{q}(-1)$ and $A_{q}(-q)$
are expressed as infinite products, \cite{Andrews1}, to properties
of and conjectures about its zeros, \cite{Andrews3,Andrews4,Hayman,Ismail5}.
It is called $q$-Airy function because it appears repeatedly in the
Plancherel-Rotach type asymptotics \cite{Ismail1,Ismail6,Ismail7}
of $q$-orthogonal polynomials, just like classical Airy function
in the classical Plancherel-Rotach asymptotics of classical orthogonal
polynomials $ $\cite{Szego,Ismail2}. In our joint work, \cite{Ismail7},
we derived Plancherel-Rotach asymptotic expansions for the $q^{-1}$-Hermite
of Ismail and Masson, $q$-Laguerre and Stieltjes-Wigert polynomials
using a discrete analogue of Laplace's method. We found that when
certain variables are above some critical values, the main terms in
the asymptotics in the bulk contain Ramanujan function $A_{q}(z)$,
when the variables are below these critical values, however, the main
terms in the asymptotics expansion in the bulk involve theta functions.

In this paper we further investigate the properties of Ramanujan function
$A_{q}(z)$. In \S\ref{sec:Preliminaries} we introduce the notations
and prove inequalities on Ismail-Masson polynomials $\left\{ h_{n}(x|q)\right\} _{n=0}^{\infty}$
and Stieltjes-Wigert polynomials $\left\{ S_{n}(x;q)\right\} _{n=0}^{\infty}$
. In \S\ref{sec:Some-Inequalities-for}, we derive two inequalities
for Ramanujan function $A_{q}(z)$. We use the asymptotic formulas
in \cite{Ismail7} to prove two definite integrals of $A_{q}(z)$
in \S\ref{sec:Definite-Integrals-for}.

\section{Preliminaries\label{sec:Preliminaries}}

In this section and the next section we will tacitly assume that all
the $\log$ and power functions are taken as their principle branches,
unless it is stated otherwise. As in our papers \cite{Ismail6,Ismail7},
we will follow the usual notations from $q$-series \cite{Andrews4,Gasper,Ismail2}

\begin{equation}
(a;q)_{0}:=1\quad(a;q)_{n}:=\prod_{k=0}^{n}(1-aq^{k}),\quad\left[\begin{array}{c}
n\\
k\end{array}\right]_{q}:=\frac{(q;q)_{n}}{(q;q)_{k}(q;q)_{n-k}}.\label{eq:2.1}\end{equation}
Thoughout this paper, we shall always assume that \begin{equation}
0<q<1\quad t>0,\label{eq:2.2}\end{equation}
hence $n=\infty$ is allowed in the above definitions. Then, \begin{equation}
0<\frac{(q;q)_{n}}{(q;q)_{n-k}}\le1\label{eq:2.3}\end{equation}
 and\begin{equation}
0<\left[\begin{array}{c}
n\\
k\end{array}\right]_{q}\le\frac{1}{(q;q)_{k}}\label{eq:2.4}\end{equation}
for $k=0,1,...,n$. 

We will use the $q$-binomial theorem \cite{Andrews4,Gasper,Ismail2},
\begin{equation}
\frac{(az;q)_{\infty}}{(z;q)_{\infty}}=\sum_{k=0}^{\infty}\frac{(a;q)_{k}}{(q;q)_{k}}z^{k}\label{eq:2.5}\end{equation}
 and the following limiting cases, also known as Euler's formulas,
\begin{equation}
(z;q)_{\infty}=\sum_{k=0}^{\infty}\frac{q^{k(k-1)/2}}{(q;q)_{k}}(-z)^{k},\qquad\frac{1}{(z;q)_{\infty}}=\sum_{k=0}^{\infty}\frac{z^{k}}{(q;q)_{k}}.\label{eq:2.6}\end{equation}
Ramanujan function $A_{q}(z)$ is defined as \cite{Ramanujan,Ismail2}
\begin{equation}
A_{q}(z):=\sum_{k=0}^{\infty}\frac{q^{k^{2}}}{(q;q)_{k}}(-z)^{k}.\label{eq:2.7}\end{equation}

\subsection{Ismail-Masson Polynomials $\left\{ h_{n}(x|q)\right\} _{n=0}^{\infty}$}

Ismail-Masson polynomials $\left\{ h_{n}(x|q)\right\} _{n=0}^{\infty}$
are defined as \cite{Ismail2} \begin{equation}
h_{n}(\sinh\xi|q)=\sum_{k=0}^{n}\left[\begin{array}{c}
n\\
k\end{array}\right]_{q}q^{k(k-n)}(-1)^{k}e^{(n-2k)\xi}.\label{eq:2.8}\end{equation}
 Ismail-Masson polynomials satisfy\begin{equation}
\int_{-\infty}^{\infty}h_{m}(x|q)h_{n}(x|q)w_{IM}(x|q)dx=q^{-n(n+1)/2}(q;q)_{n}\delta_{m,n}\label{eq:2.9}\end{equation}
 for $n,m=0,1,...$, where \begin{equation}
w_{IM}(x|q)=\sqrt{\frac{-2q^{1/4}}{\pi\log q}}\,\exp\left\{ \frac{2}{\log q}\left[\log\left(x+\sqrt{x^{2}+1}\right)\right]^{2}\right\} .\label{eq:2.10}\end{equation}
It is clear that their orthonormal polynomials are defined as \begin{equation}
\tilde{h}_{n}(x|q)=\frac{q^{n(n+1)/4}}{\sqrt{(q;q)_{n}}}h_{n}(x|q).\label{eq:2.11}\end{equation}
Let \begin{equation}
\sinh\xi_{n}:=\frac{q^{-nt}u-q^{nt}u^{-1}}{2},\label{eq:2.12}\end{equation}
and assume that \begin{equation}
u\in\mathbb{C}\backslash\left\{ 0\right\} .\label{eq:2.13}\end{equation}
 It is easy to see that\begin{equation}
w_{IM}(\sinh\xi_{n}|q)=w_{IM}(\sinh u|q)u^{-4nt}q^{2n^{2}t^{2}}.\label{eq:2.14}\end{equation}
It is also clear from \eqref{eq:2.8} and \eqref{eq:2.12} that

\begin{equation}
h_{n}(\sinh\xi_{n}|q)=u^{n}q^{-n^{2}t}\sum_{k=0}^{n}\left[\begin{array}{c}
n\\
k\end{array}\right]_{q}q^{k^{2}}\left(-\frac{q^{n(2t-1)}}{u^{2}}\right)^{k}.\label{eq:2.15}\end{equation}
Thus

\begin{eqnarray*}
\left|h_{n}(\sinh\xi_{n}|q)\right| & \le & \frac{\left|u\right|^{n}}{q^{n^{2}t}}\sum_{k=0}^{n}\frac{q^{k^{2}}}{(q;q)_{k}}\left(\frac{q^{n(2t-1)}}{\left|u\right|^{2}}\right)^{k}\\
 & \le & \frac{\left|u\right|^{n}}{q^{n^{2}t}}\sum_{k=0}^{\infty}\frac{q^{k^{2}}}{(q;q)_{k}}\left(\frac{q^{n(2t-1)}}{\left|u\right|^{2}}\right)^{k},\end{eqnarray*}
or\begin{equation}
\left|h_{n}(\sinh\xi_{n}|q)\right|\le\frac{\left|u\right|^{n}}{q^{n^{2}t}}A_{q}\left(-\frac{q^{n(2t-1)}}{\left|u\right|^{2}}\right).\label{eq:2.16}\end{equation}

\subsection{Stieltjes-Wigert Polynomials $\left\{ S_{n}(x;q)\right\} _{n=0}^{\infty}$ }

Stieltjes-Wigert polynomials $\left\{ S_{n}(x;q)\right\} _{n=0}^{\infty}$
are defined as \cite{Ismail2}\begin{equation}
S_{n}(x;q)=\sum_{k=0}^{n}\frac{q^{k^{2}}(-x)^{k}}{(q;q)_{k}(q;q)_{n-k}}.\label{eq:2.17}\end{equation}
They are orthogonal respect to the weight function\begin{equation}
w_{SW}(x;q)=\sqrt{\frac{-1}{2\pi\log q}}\exp\left\{ \frac{1}{2\log q}\left[\log\left(\frac{x}{\sqrt{q}}\right)\right]^{2}\right\} ,\label{eq:2.18}\end{equation}
with\begin{equation}
\int_{0}^{\infty}S_{n}(x;q)S_{m}(x;q)w_{SW}(x;q)dx=\frac{q^{-n}}{(q;q)_{n}}\delta_{m,n}\label{eq:2.19}\end{equation}
for $n,m=0,1,...$. The orthonormal Stieltjes-Wigert polynomials with
positive leading coefficients are \begin{equation}
\tilde{s}_{n}(x;q)=(-1)^{n}\sqrt{q^{n}(q;q)_{n}}\, S_{n}(x;q).\label{eq:2.20}\end{equation}
In the case of the Stieltjes-Wigert polynomials the appropriate scaling
is \begin{equation}
x_{n}(t,u)=q^{-nt}u.\label{eq:2.21}\end{equation}
 A calculation gives\begin{equation}
w_{SW}(q^{-nt}u;q)=w_{SW}(u;q)u^{-nt}q^{(n^{2}t^{2}+nt)/2}.\label{eq:2.22}\end{equation}
Set $x=x_{n}(t,u)$ in \eqref{eq:2.17} then replace $k$ by $n-k$
to see that \begin{equation}
S_{n}(x_{n}(t,u);q)=u^{n}q^{n^{2}(1-t)}\sum_{k=0}^{n}\frac{q^{k^{2}}\left(-\frac{q^{n(t-2)}}{u}\right)^{k}}{(q;q)_{k}(q;q)_{n-k}}.\label{eq:2.23}\end{equation}
Thus,\begin{eqnarray*}
\left|S_{n}(x_{n}(t,u);q)(q;q)_{n}\right| & \le & \frac{\left|u\right|^{n}}{q^{n^{2}(t-1)}}\sum_{k=0}^{n}\frac{q^{k^{2}}}{(q;q)_{k}}\left(\frac{q^{n(t-2)}}{\left|u\right|}\right)^{k}\\
 & \le & \frac{\left|u\right|^{n}}{q^{n^{2}(t-1)}}\sum_{k=0}^{\infty}\frac{q^{k^{2}}}{(q;q)_{k}}\left(\frac{q^{n(t-2)}}{\left|u\right|}\right)^{k},\end{eqnarray*}
 or \begin{equation}
\left|S_{n}(x_{n}(t,u);q)\right|\le\frac{\left|u\right|^{n}A_{q}\left(-\frac{q^{n(t-2)}}{\left|u\right|}\right)}{(q;q)_{\infty}q^{n^{2}(t-1)}}.\label{eq:2.24}\end{equation}

\section{Some Inequalities for $A_{q}(z)$\label{sec:Some-Inequalities-for}}

It is clear that\begin{equation}
n\ge\frac{1-q^{n}}{1-q}\ge nq^{n-1}\label{eq:3.1}\end{equation}
for $n\in\mathbb{N}$, then, \begin{equation}
\left|\frac{(1-q)^{k}}{(q;q)_{k}}q^{k^{2}}(-z)^{k}\right|\le\frac{(q\left|z\right|)^{k}}{k!}\label{eq:3.2}\end{equation}
 for $k=0,1,\dotsc$ and for any complex number $z$. Applying Lebesgue's
dominated convergent theorem we have\begin{equation}
\lim_{q\to1}A_{q}((1-q)z)=e^{-z}\label{eq:3.3}\end{equation}
 for any $z\in\mathbb{C}$, hence $A_{q}(z)$ is really one of many
$q$-analogues of the exponential function. From \eqref{eq:3.2} we
also have obtained the inequality \begin{equation}
\left|A_{q}((1-q)z)\right|\le e^{q|z|}\label{eq:3.4}\end{equation}
 or\begin{equation}
\left|A_{q}(z)\right|\le e^{q|z|/(1-q)}\label{eq:3.5}\end{equation}
for any complex number $z$. For any nonzero complex number $z$,
then\begin{eqnarray}
\left|A_{q}(z)\right| & \le & \sum_{k=0}^{\infty}\frac{q^{k}}{(q;q)_{k}}\left(q^{k-1}\left|z\right|\right)^{k}.\label{eq:3.6}\end{eqnarray}
For $k=0,1,...$, the terms $q^{k(k-1)}|z|^{k}$ are bounded by \begin{equation}
\left(\frac{|z|}{\sqrt{q}}\right)^{1/2}\exp\left\{ -\frac{\log^{2}\left|z\right|}{4\log q}\right\} .\label{eq:3.7}\end{equation}
 Thus,

\begin{eqnarray*}
\left|A_{q}(z)\right| & \le & \left(\frac{|z|}{\sqrt{q}}\right)^{1/2}\exp\left\{ -\frac{\log^{2}\left|z\right|}{4\log q}\right\} \sum_{k=0}^{\infty}\frac{q^{k}}{(q;q)_{k}}\\
 & \le & \frac{\left(\frac{|z|}{\sqrt{q}}\right)^{1/2}\exp\left\{ -\frac{\log^{2}\left|z\right|}{4\log q}\right\} }{(q;q)_{\infty}}\end{eqnarray*}
 or we have\begin{equation}
\left|A_{q}(z)\right|\le\frac{\left(\frac{|z|}{\sqrt{q}}\right)^{1/2}\exp\left\{ -\frac{\log^{2}\left|z\right|}{4\log q}\right\} }{(q;q)_{\infty}}\label{eq:3.8}\end{equation}
 for any nonzero complex number $z$. 

\begin{thm}
Assume that $A_{q}(z)$ is Ramanujan function defined in \eqref{eq:2.6},
then, for any complex number $z$, \begin{equation}
\left|A_{q}(z)\right|\le e^{q|z|/(1-q)},\label{eq:3.9}\end{equation}
 and , \begin{equation}
\left|A_{q}(z)\right|\le\frac{\left(\frac{|z|}{\sqrt{q}}\right)^{1/2}\exp\left\{ -\frac{\log^{2}\left|z\right|}{4\log q}\right\} }{(q;q)_{\infty}},\label{eq:3.10}\end{equation}
 for any complex number $z\neq0$. 
\end{thm}
\begin{rem}
The trivial inequality \eqref{eq:3.1} can be used to show that a
basic hypergeometric series converges to its hypergeometric series
counter-part under suitable scaling and conditions. Also, using \eqref{eq:3.10}
the formulas \eqref{eq:2.16} and \eqref{eq:2.24} could be recast
into other forms.
\end{rem}

\section{Definite Integrals for $A_{q}(z)$\label{sec:Definite-Integrals-for}}

Put $t=\frac{1}{2}$ in formula (64) of \cite{Ismail7}, we have\begin{eqnarray}
\sqrt{\frac{(q;q)_{n}w_{H}(\sinh\xi_{n}|q)}{q^{n/2}w_{H}(\sinh u|q)}}\tilde{h}_{n}(\sinh\xi_{n}|q) & = & A_{q}\left(u^{-2}\right)+r_{IM}\label{eq:4.1}\end{eqnarray}
 with \begin{gather}
\left|r_{IM}\right|\le\frac{4(-q^{3};q)_{\infty}A_{q}\left(-|u|^{-2}\right)}{(q;q)_{\infty}^{2}}\label{eq:4.2}\\
\times\left(q^{n/2}+q^{n^{2}/4}|u|^{-2\lfloor n/2\rfloor-2}\right).\nonumber \end{gather}
 \[
\]

\begin{thm}
Assuming that $A_{q}(z)$ and $w_{IM}(x|q)$ are defined as in \eqref{eq:2.7}
and \eqref{eq:2.10}. Then,\begin{equation}
\int_{0}^{\infty}A_{q}^{2}(u^{-2})w_{IM}(u|q)du=2(q;q)_{\infty}.\label{eq:4.3}\end{equation}
 
\end{thm}
\begin{proof}
For the orthonormal Ismail-Masson polynomials $\tilde{h}_{n}(x|q)$
$ $satisfy\begin{equation}
\int_{-\infty}^{\infty}\left\{ \tilde{h}_{n}(x|q)\right\} ^{2}w_{IM}(x|q)dx=1.\label{eq:4.4}\end{equation}
 Assume $u>0$ and make the change of variable \begin{equation}
x=\sinh\xi_{n}=\frac{q^{-n/2}u-q^{n/2}u^{-1}}{2}\label{eq:4.5}\end{equation}
in \eqref{eq:4.4}, we have\begin{equation}
\int_{0}^{\infty}\left\{ \tilde{h}_{n}(\sinh\xi_{n}|q)\right\} ^{2}w_{IM}(\sinh\xi_{n}|q)(q^{-n/2}+q^{n/2}u^{-2})du=2,\label{eq:4.6}\end{equation}
 or\begin{equation}
\int_{0}^{\infty}\left\{ \sqrt{\frac{(q;q)_{n}w_{IM}(\sinh\xi_{n}|q)}{q^{n/2}w_{IM}(\sinh u|q)}}\,\tilde{h}_{n}(\sinh\xi_{n}|q)\right\} ^{2}\left(1+q^{n}u^{-2}\right)w_{IM}(\sinh u|q)du=2(q;q)_{n}.\label{eq:4.7}\end{equation}
Thus we have\begin{equation}
\lim_{n\to\infty}\int_{0}^{\infty}\left\{ \sqrt{\frac{(q;q)_{n}w_{IM}(\sinh\xi_{n}|q)}{q^{n/2}w_{IM}(\sinh u|q)}}\,\tilde{h}_{n}(\sinh\xi_{n}|q)\right\} ^{2}\left(1+q^{n}u^{-2}\right)w_{IM}(\sinh u|q)du=2(q;q)_{\infty}.\label{eq:4.8}\end{equation}
From \eqref{eq:2.11}, \eqref{eq:2.14} and \eqref{eq:2.16} we have
\begin{equation}
\left\{ \sqrt{\frac{(q;q)_{n}w_{IM}(\sinh\xi_{n}|q)}{q^{n/2}w_{IM}(\sinh u|q)}}\,\left|\tilde{h}_{n}(\sinh\xi_{n}|q)\right|\right\} ^{2}\le A_{q}^{2}(-u^{-2}),\label{eq:4.9}\end{equation}
 and from \eqref{eq:3.10} and \eqref{eq:2.10} we know that \begin{equation}
u^{-2}w_{IM}(\sinh u|q)A_{q}^{2}(-u^{-2})\label{eq:4.10}\end{equation}
is bounded for $0<u\le1$. From \eqref{eq:3.9} we know that\begin{equation}
A_{q}^{2}(-u^{-2})\label{eq:4.11}\end{equation}
 is bounded for $u\ge1$. Lebesgue dominated convergence theorem allows
us to take limit inside the integral. We use \eqref{eq:4.1} to get
\eqref{eq:4.3}.
\end{proof}
Take $t=2$ in the formula (69) of \cite{Ismail7}, 

\begin{eqnarray}
\sqrt{\frac{q^{-n}w_{SW}(q^{-2n}u;q)}{(q;q)_{n}w_{SW}(u;q)}}\tilde{s}_{n}(q^{-2n}u;q) & = & \frac{\left\{ A_{q}(u^{-1})+r_{SW}(n)\right\} }{(q;q)_{\infty}}\label{eq:4.12}\end{eqnarray}
with\begin{equation}
|r_{SW}(n)|\le\frac{2(-q^{3};q)_{\infty}A_{q}\left(-|u|^{-1}\right)}{(q;q)_{\infty}}\left\{ q^{n/2}+\frac{q^{n^{2}/4}}{|u|^{1+\lfloor n/2\rfloor}}\right\} .\label{eq:4.13}\end{equation}

\begin{thm}
Assuming that $A_{q}(z)$, and $w_{SW}(x;q)$ are defined as in \eqref{eq:2.7}
and \eqref{eq:2.18}. Then we have\begin{equation}
\int_{0}^{\infty}A_{q}^{2}(u^{-1})w_{SW}(u;q)du=(q;q)_{\infty}.\label{eq:4.14}\end{equation}

\end{thm}
\begin{proof}
From the orthogonality of Stieltjes-Wigert polynomials, we know that
the orthonormal polynomials $\tilde{s}_{N}(x;q)$ satisfy\begin{equation}
\int_{0}^{\infty}\tilde{s}_{n}^{2}(x;q)w_{SW}(x;q)dx=1.\label{eq:4.15}\end{equation}
Let us make a change of variable \begin{equation}
x=q^{-2n}u\label{eq:4.16}\end{equation}
 in \eqref{eq:4.15} with $u>0$, then,\begin{equation}
\int_{0}^{\infty}\left[\sqrt{\frac{q^{-2n}w_{SW}(q^{-2n}u;q)}{(q;q)_{n}w_{SW}(u;q)}}\,\tilde{s}_{n}(q^{-2n}u;q)\right]^{2}w_{SW}(u;q)du=\frac{1}{(q;q)_{n}},\label{eq:4.17}\end{equation}
 therefore,\begin{equation}
\lim_{n\to\infty}\int_{0}^{\infty}\left[\sqrt{\frac{q^{-2n}w_{SW}(q^{-2n}u;q)}{(q;q)_{n}w_{SW}(u;q)}}\,\tilde{s}_{n}(q^{-2n}u;q)\right]^{2}w_{SW}(u;q)du=\frac{1}{(q;q)_{\infty}}.\label{eq:4.18}\end{equation}
 From \eqref{eq:2.20}, \eqref{eq:2.22} and \eqref{eq:2.24} we have\[
\left\{ \sqrt{\frac{q^{-2n}w_{SW}(q^{-2n}u;q)}{(q;q)_{n}w_{SW}(u;q)}}\,\left|\tilde{s}_{n}(q^{-2n}u;q)\right|\right\} ^{2}\le\frac{A_{q}^{2}(-u^{-1})}{(q;q)_{\infty}^{2}},\]
 and by \eqref{eq:2.14} and \eqref{eq:3.10}\[
A_{q}^{2}(-u^{-1})w_{SW}(u;q)\]
 is bounded for $0<u\le1$, by \eqref{eq:3.9}, \[
A_{q}^{2}(-u^{-2})\]
 is bounded for $u\ge1$. By Lebesgue dominated convergence theorem
we interchange the orders of limit and integration, then apply \eqref{eq:4.12}
to get \eqref{eq:4.14}.
\end{proof}


\begin{thebibliography}{10}
\bibitem{Akhiezer}N. I. Akhiezer, The Classical Moment Problem and
Some Related Questions in Analysis, English translation, Oliver and
Boyed, Edinburgh, 1965.

\bibitem{Andrews1}G. E. Andrews, q-series: Their development and
application in analysis, number theory, combinatorics, physics, and
computer algebra, CBMS Regional Conference Series, number 66, American
Mathematical Society, Providence, R.I. 1986.

\bibitem{Andrews2}G. E. Andrews, Ramanujan's \char`\"{}Lost\char`\"{}
Note book VIII: The entire Rogers-Ramanujan function, Advances in
Math. 191 (2005), 393--407.

\bibitem{Andrews3}G. E. Andrews, Ramanujan's \char`\"{}Lost\char`\"{}
Note book IX: The entire Rogers-Ramanujan function, Advances in Math.
191 (2005), 408--422.

\bibitem{Andrews4}G. E. Andrews, R. A. Askey, and R. Roy, Special
Functions, Cambridge University Press, Cambridge, 1999.

\bibitem{Deift1}P. Deift, Orthogonal Polynomials and Random Matrices:
a Riemann-Hilbert Approach, American Mathematical Society, Providence,
2000.

\bibitem{Deift2}P. Deift, T. Kriecherbauer, K. T-R. McLaughlin, S.
Venakides, and X. Zhou, Strong asymptotics of orthogonal polynomials
with respect to exponential weights, Comm. Pure Appl. Math. 52 (1999),
1491--1552.

\bibitem{Gasper}G. Gasper and M. Rahman, Basic Hypergeometric Series,
second edition Cambridge University Press, Cambridge, 2004.

\bibitem{Hayman}W. K. Hayman, On the zeros of a q-Bessel function,
Contemporary Mathematics, volume 382, American Mathematical Society,
Providence, 2005, 205--216.

\bibitem{Ismail1}M. E. H. Ismail, Asymptotics of q-orthogonal polynomials
and a q-Airy function, Internat. Math. Res. Notices 2005 No 18 (2005),
1063--1088.

\bibitem{Ismail2}M. E. H. Ismail, Classical and Quantum Orthogonal
Polynomials in one Variable, Cambridge University Press, Cambridge,
2005.

\bibitem{Ismail3}M. E. H. Ismail and X. Li, Bounds for extreme zeros
of orthogonal polynomials, Proc. Amer. Math. Soc. 115 (1992), 131--140.

\bibitem{Ismail4}M. E. H. Ismail and D. R. Masson, q-Hermite polynomials,
biorthogonal rational functions, Trans. Amer. Math. Soc. 346 (1994),
63--116.

\bibitem{Ismail5}M. E. H. Ismail and C. Zhang, Zeros of entire functions
and a problem of Ramanujan, Advances in Math., (2007), to appear.

\bibitem{Ismail6}M. E. H. Ismail and R. Zhang, Scaled asymptotics
for q-polynomials, Comptes Rendus, submitted.

\bibitem{Ismail7}M. E. H. Ismail and R. Zhang, Chaotic and Periodic
Asymptotics for q-Orthogonal Polynomials, joint with Mourad E.H. Ismail,
International Mathematics Research Notices, accepted. 

\bibitem{Kajiwara}K. Kajiwara, T. Masuda, M. Noumi, Y. Ohta, Y. Yamada,
Hypergeometric solutions to the q-Painlev\textbackslash{}'\{e\} equations,
Internat. Math. Res. Notices 47 (2004), 2497--2521.

\bibitem{Koekoek}R. Koekoek and R. Swarttouw, The Askey-scheme of
hypergeometric orthogonal polynomials and its q-analogues, Reports
of the Faculty of Technical Mathematics and Informatics no. 98-17,
Delft University of Technology, Delft, 1998.

\bibitem{Mehta}M. L. Mehta, Random Matrices, third edition, Elsevier,
Amsterdam, 2004.

\bibitem{Ramanujan}S. Ramanujan, The Lost Notebook and Other Unpublished
Papers (Introduction by G. E. Andrews), Narosa, New Delhi, 1988.

\bibitem{Saff}E. B. Saff and V. Totik, Logarithmic Potentials With
External Fields, Springer-Verlag, New York, 1997.

\bibitem{Szego}G. Szeg\textbackslash{}H\{o\}, Orthogonal Polynomials,
Fourth Edition, Amer. Math. Soc., Providence, 1975.

\bibitem{Wong}R. Wong, Asymptotic Approximations of Integrals, Academic
Press, Boston, 1989.

\bibitem{Whattaker}E. T. Whittaker and G. N. Watson, A Course of
Modern Analysis, fourth edition, Cambridge University Press, Cambridge,
1927. 
\end{thebibliography}
\end{document}